\documentclass[oneside]{amsart}

\newcommand{\hide}[1]{}

\usepackage{amssymb}
\usepackage{graphicx}
\usepackage{amsmath}

\usepackage{array} 

\newcolumntype{P}[1]{>{\centering\arraybackslash}p{#1}}
\newcolumntype{C}[1]{>{\centering\arraybackslash}S{m{#1}}}

\usepackage{cellspace}
\usepackage{mathrsfs}
\usepackage{geometry}

\usepackage{tikz-cd}

\usepackage{xcolor}

\usepackage{amsfonts,subcaption}
\usepackage{amssymb}
\usepackage{enumitem}

\usepackage[colorlinks=true,urlcolor=blue]{hyperref}

\DeclareFontFamily{U}{mathb}{\hyphenchar\font45}
\DeclareFontShape{U}{mathb}{m}{n}{ <5> <6> <7> <8> <9> <10> gen * mathb <10.95> mathb10 <12> <14.4> <17.28> <20.74> <24.88> mathb12 }{}
\DeclareSymbolFont{mathb}{U}{mathb}{m}{n}
\DeclareFontSubstitution{U}{mathb}{m}{n}
\DeclareMathSymbol{\righttoleftarrow}{3}{mathb}{"FD}
\DeclareMathSymbol{\selfmap}{3}{mathb}{"FD}

\makeatletter
\newcommand{\setword}[2]{%
  \phantomsection
  #1\def\@currentlabel{\unexpanded{#1}}\label{#2}%
}

\makeatletter
\renewcommand*\env@matrix[1][*\c@MaxMatrixCols c]{%
  \hskip -\arraycolsep
  \let\@ifnextchar\new@ifnextchar
  \array{#1}}
\makeatother

\newcommand{\R}{\mathbb{R}}
\newcommand{\C}{\mathbb{C}}

\renewcommand{\H}{\mathbb{H}}

\newcommand{\Cc}{\widehat{{\C}}}

\DeclareMathOperator{\diver}{div}

\newcommand{\cmod}[1]{\modulus\left(#1\right)}

\DeclareMathOperator{\Area}{Area}

\renewcommand{\rho}{\varrho}
\renewcommand{\phi}{\varphi}
\newcommand{\eps}{\varepsilon}

\renewcommand{\theta}{\vartheta}

\DeclareMathOperator{\modulus}{mod}

\DeclareMathOperator{\dist}{dist}

\newcommand\ovl[1]{\overline{#1}}
\newcommand{\sm}{\setminus}

\renewcommand{\ge}{\geqslant}
\renewcommand{\le}{\leqslant}

\theoremstyle{plain}
\newtheorem{theorem}{Theorem}
\newtheorem{lemma}[theorem]{Lemma}

\theoremstyle{definition}

\newtheorem{remark}{Remark}

\everymath{\displaystyle}
\definecolor{fgreen}{rgb}{0.2,0.6,0.}

\newcounter{reminder}

\newtheoremstyle{claim}
  {}
  {}
  {\itshape}
  {0pt}
  {\scshape}
  {.}
  { }
  {\thmname{#1}\thmnumber{ #2}\thmnote{ (#3)}}

\theoremstyle{claim}

\numberwithin{equation}{section}

\title[Hyperbolic area proof of the Kahn--Lyubich covering lemma]{A hyperbolic area proof of the Kahn--Lyubich Covering Lemma and its generalization}
\author{Kostiantyn Drach}

\date{}

\thanks{This work was partially supported by Agencia Estatal de Investigaci\'on grants PID2023-147252NB-I00, CNS2025-166633, and by the Ag\`encia de Gesti\'o d'Ajuts Universitaris i de Recerca of the Generalitat de Catalunya grant number 2025 SGR 00322.}

\address{Departament de Matem{\`a}tiques i Inform{\`a}tica, Universitat de Barcelona, Gran Via de les Corts Catalanes 585, 08007 Barcelona, Spain}
\address{Centre de Recerca Matem{\`a}tica, Edifici C, Carrer de l'Albareda, 08193 Bellaterra, Spain}

\email{kostiantyn.drach@ub.edu}

\begin{document}

\begin{abstract}
    The Covering Lemma of J.~Kahn and M.~Lyubich is a fundamental analytic tool for controlling conformal moduli of annuli under branched coverings of finite degree. In the last several decades, it has played a crucial role in establishing rigidity of multicritical polynomial and rational maps. The classical proof uses the Near-Degenerate Regime as the main mechanism. In this short note, we prove this lemma using hyperbolic area and the Dirichlet energy. We also establish a version of the Covering Lemma for holomorphic maps of infinite degree. 
    
    Our argument is inspired by recent applications of the hyperbolic area method to the study of wandering domains in rational and transcendental dynamics.  
\end{abstract}

\maketitle

\section{Introduction}

The Covering Lemma of J.~Kahn and M.~Lyubich \cite{KL} has been essential in establishing local connectivity and rigidity results for complex polynomials of higher degree \cite{KLLC, KvS, DY} and, through puzzle renormalization via complex box mappings \cite{KvS, CDKS}, for rational maps, most notably Newton maps \cite{DS, WYZ, RYZ}. It has also been used to prove local connectivity of boundaries of bounded polynomial Fatou components \cite{RY, Yang}. Furthermore, the Covering Lemma has played an important role in recent advances towards one of the central conjectures in the field -- the MLC conjecture (``the Mandelbrot set is locally connected"), and in the development of the Near-Degenerate Regime mechanism; we refer the reader to a very nice survey \cite{DudkoICM}.

The original proof of the Covering Lemma \cite{KL} relies on the near-degenerate philosophy and a careful analysis of conformal rectangles and their moduli (extremal widths). Inspired by recent applications of the hyperbolic area method to wandering domains in rational and transcendental dynamics \cite{DPU,PRW, ye26-2, ye26}, in this note we present an alternative (somewhat elementary) proof using hyperbolic geometry, hyperbolic area, and its connection with conformal moduli and extremal widths through Dirichlet energy. Furthermore, we extend the Covering Lemma to maps of infinite degree.

We consider non-constant holomorphic maps $f \colon U \to V$ between topological disks with finitely many \emph{singular values}. Thus, there is a finite set $S_f \subset V$ such that $f \colon U \sm f^{-1}(S_f) \to V \sm S_f$ is a covering map; we choose $S_f$ minimal with this property. The elements of $S_f$ are \emph{critical} or \emph{asymptotic} values of $f$. If $f$ is proper, then it has finite degree and its singular values are precisely its critical values.

\newpage
\begin{theorem}[The Covering Lemma]
    \label{Thm:Main}
    Let $\Lambda \Subset \Lambda' \Subset U$ and $B \Subset B' \Subset V$ be two nests of open Jordan disks in $\Cc$. Let $f \colon U \to V$ be a non-constant holomorphic map with at most $D \ge 1$ singular values. Assume that the restrictions $f \colon \Lambda' \to B'$ and $f \colon \Lambda \to B$ are proper, and that the degree $\deg(\Lambda' \to B') = d \le D$.

    Then for every $\eta \in (0,1]$ there exists an explicit $\eps = \eps(\eta,D)$ with the following property: If the annulus $U \sm \ovl \Lambda$ is $\eps$-degenerate, i.e.,
    \[
    \modulus\left(U \sm \ovl \Lambda\right) < \eps,
    \]
    then one of the following alternatives holds:
    \begin{enumerate}
        \item either 
        \[
        \modulus\left(U \sm \ovl \Lambda\right) > \frac{1}{\eta} \modulus\left(B' \sm \ovl B\right),
        \]
        \item or the \emph{collar assumption} $\modulus\left(B' \sm \ovl B\right) \ge \eta \modulus\left(U \sm \ovl \Lambda\right)$ is satisfied and
        \[
        \modulus\left(U \sm \ovl \Lambda\right) > \frac{\eta}{32d^2} \modulus\left(V \sm \ovl B\right).
        \]
    \end{enumerate}
\end{theorem}

\begin{figure}[t]
\begin{center}
\definecolor{grey1}{rgb}{0.5,0.5,0.5}
\definecolor{grey2}{rgb}{0.25,0.25,0.25}
\begin{tikzpicture}[line cap=round,line join=round,>=stealth,x=1.0cm,y=1.0cm]
\clip(-4.4,-2.7) rectangle (8.7,3.);
\draw[line width=1.pt,color=grey1,fill=grey1,fill opacity=0.1] (-4.239886417383266,0.4956662291311774) -- (-4.239886415498974,0.495666422739675) -- (-4.239886411730362,0.4956668099566538) -- (-4.239886404193024,0.49566758439054553) -- (-4.239886389117893,0.4956691332580658) -- (-4.239886358965814,0.49567223099205343) -- (-4.239886298654389,0.4956784264558169) -- (-4.239886178002463,0.495690817366497) -- (-4.239885936582316,0.49571559912046875) -- (-4.239885453276842,0.49576516235885276) -- (-4.239884484805175,0.4958642877573356) -- (-4.2398825404189795,0.49606253424078434) -- (-4.239878621875287,0.4964590099506033) -- (-4.239870665703816,0.49725189231784245) -- (-4.23985427703348,0.4988373806500432) -- (-4.239819594454915,0.5020072501638797) -- (-4.239742608919623,0.5083425482576565) -- (-4.239558160923663,0.5209952820565281) -- (-4.239067393928455,0.5462285108747752) -- (-4.237598669510539,0.596399699595128) -- (-4.23271480765,0.6955104909860222) -- (-4.227391547381027,0.7682301044043103) -- (-4.220632355797756,0.8399370036522291) -- (-4.212439453228114,0.910583438773884) -- (-4.20281506000003,0.9801216598133801) -- (-4.179280682880243,1.1156824598223158) -- (-4.150046987061813,1.246237404031878) -- (-4.11513173516816,1.3714044927949076) -- (-4.074552689822701,1.4908017264642455) -- (-4.028327613648856,1.6040471053927328) -- (-3.976474269270044,1.710758629933211) -- (-3.9190104193096844,1.8105543004385205) -- (-3.855953826391196,1.9030521172615031) -- (-3.787322253137997,1.9878700807549994) -- (-3.713133462173507,2.0646261912718513) -- (-3.6334052161211448,2.1329384491648993) -- (-3.5481557356591105,2.1924257647462104) -- (-3.457500235612125,2.2428997342752455) -- (-3.361770281860552,2.2846017511653365) -- (-3.261326705370645,2.317831346057707) -- (-3.1565303371086557,2.342888049593578) -- (-3.0477420080408377,2.3600713924141736) -- (-2.935322549133444,2.369680905160715) -- (-2.819632791352728,2.3720161184744257) -- (-2.7010335656649427,2.3673765629965278) -- (-2.5798857030363394,2.3560617693682424) -- (-2.4565500344331728,2.3383712682307944) -- (-2.3313873908216953,2.3146045902254047) -- (-2.2047586031681603,2.2850612659932965) -- (-2.07702450243882,2.2500408261756917) -- (-1.9485459195999288,2.2098428014138136) -- (-1.819683685617738,2.164766722348883) -- (-1.690798631458502,2.1151121196221236) -- (-1.5622515880884724,2.0611785238747595) -- (-1.4344033864739032,2.00326546574801) -- (-1.3076148575810462,1.9416724758830979) -- (-1.1822468323761566,1.8766990849212464) -- (-1.0586601418254857,1.8086448235036814) -- (-0.9372143015003376,1.7378070825717296) -- (-0.8182046886628331,1.6643789218660898) -- (-0.7018349469711254,1.5884041817739307) -- (-0.588301663122488,1.5099152234087836) -- (-0.4778014238141939,1.4289444078841786) -- (-0.3705308157435194,1.3455240963136463) -- (-0.26668642560774014,1.259686649810717) -- (-0.16646484010412976,1.1714644294889225) -- (-0.07006264592996381,1.080889796461792) -- (0.022323570217485766,0.9879951118428565) -- (0.11049722164094256,0.8928127367456455) -- (0.19426172164313016,0.7953750322836908) -- (0.2734204835267766,0.6957143595705224) -- (0.3477769205946055,0.5938630797196711) -- (0.41713444614934403,0.48985355384466756) -- (0.48129647349371485,0.3837181430590415) -- (0.5400664159304442,0.2754892084763241) -- (0.5932476867622603,0.1651991112100446) -- (0.6406436992918865,0.05288021237373508) -- (0.6820578668220456,-0.06143512691907338) -- (0.7172936026554666,-0.17771454555485233) -- (0.7461543200948713,-0.29592568242007067) -- (0.7684434324429885,-0.4160361764011964) -- (0.7839643530025411,-0.5380136663847006) -- (0.7925204950762579,-0.6618257912570513) -- (0.7939152719668607,-0.7874401899047196) -- (0.7879520969770741,-0.9148245012141745) -- (0.774434383409627,-1.0439463640718847) -- (0.7531655445672438,-1.1747734173643218) -- (0.7239565939106143,-1.3072571460184434) -- (0.6869146945227129,-1.4407195758042661) -- (0.642532119991202,-1.5736641890422192) -- (0.5913269219311843,-1.7045396774526722) -- (0.5338171519579333,-1.8317947327557675) -- (0.47052086168656615,-1.9538780466721022) -- (0.40195610273234195,-2.06923831092179) -- (0.3286409267103636,-2.1763242172253143) -- (0.2510933852358619,-2.2735844573027606) -- (0.16983152992402495,-2.359467722874726) -- (0.08537341239005514,-2.432422705661324) -- (-0.0017629157509162496,-2.4908980973830666) -- (-0.09105940288365844,-2.5333425897601245) -- (-0.18199902868099826,-2.5582103649494456) -- (-0.2742480079987786,-2.5649311242741533) -- (-0.36786522701536484,-2.5550250977112228) -- (-0.4629602867411684,-2.530282514098417) -- (-0.5596427881865953,-2.492493602273953) -- (-0.6580223323620533,-2.4434485910753097) -- (-0.758208520277952,-2.3849377093406474) -- (-0.8603109529446993,-2.3187511859075585) -- (-0.9644392313727028,-2.2466792496140897) -- (-1.0707029565723705,-2.170512129298004) -- (-1.1792117295541065,-2.092040053797291) -- (-1.2900751513283257,-2.0130532519497137) -- (-1.4034028229054432,-1.9353419525932622) -- (-1.5193043452958488,-1.8606963845655855) -- (-1.6378893195099575,-1.7909067767047873) -- (-1.7592673465581878,-1.7277633578485165) -- (-1.8835480274509404,-1.6730563568348202) -- (-2.0108379950274013,-1.628565299403931) -- (-2.1409420870760414,-1.5949814512964622) -- (-2.2731075743047313,-1.5709855152966838) -- (-2.406521480787262,-1.5550409474109301) -- (-2.5403708305973964,-1.545611203646331) -- (-2.673842647808897,-1.5411597400092205) -- (-2.8061239564953553,-1.5401500125063876) -- (-2.9364017807305345,-1.54104547714428) -- (-3.0638631445883675,-1.5423095899298005) -- (-3.1876950721423896,-1.5424058068698514) -- (-3.30708458746642,-1.5397975839705396) -- (-3.421218714634108,-1.5329483772386538) -- (-3.5292844777193295,-1.5203216426813242) -- (-3.6304689007957336,-1.5003808363045437) -- (-3.7239590079369123,-1.4715894141154422) -- (-3.8089418232168555,-1.4324108321208087) -- (-3.884661129426746,-1.3814016911775866) -- (-3.9511353618106284,-1.3183898486219618) -- (-4.008943282743161,-1.244122696080069) -- (-4.058676133612295,-1.1593681040739057) -- (-4.1009251558063795,-1.0648939431260374) -- (-4.136281590713139,-0.961468083758291) -- (-4.165336679721037,-0.8498583964932891) -- (-4.18868166421791,-0.7308327518528586) -- (-4.206907785592051,-0.6051590203592241) -- (-4.220606285231298,-0.4736050725348946) -- (-4.230368404523944,-0.3369387789016969) -- (-4.23395808447404,-0.26692846194498543) -- (-4.2367853848579955,-0.19592800998196935) -- (-4.238924210849433,-0.1240334068280049) -- (-4.240448467621633,-0.05134063629799357) -- (-4.241432060348217,0.02205431779253786) -- (-4.241948894202807,0.0960554716280626) -- (-4.242072874358797,0.170566841393736) -- (-4.241877905989696,0.2454924432740313) -- (-4.2414378942692395,0.32073629345387644) -- (-4.240826744370537,0.39620240811774465) -- (-4.240118361467495,0.4717948034505639) -- (-4.239934733847264,0.49070138553503284) -- (-4.239888734455718,0.49542815117223427) -- (-4.239887296866215,0.49557586273925835) -- (-4.239886578071292,0.4956497185229978) -- (-4.23988648822197,0.4956589504959652) -- (-4.2398864432973085,0.49566356648239207) -- (-4.2398864208348925,0.4956658744757192)(-4.105761435343545,1.4010555035331393);
\draw[line width=1.pt,color=fgreen,fill=fgreen,fill opacity=0.1] (-2.212179672707263,0.5598787160656278) -- (-2.212179646510154,0.5598787974261312) -- (-2.212179594115908,0.5598789601471196) -- (-2.2121794893273,0.5598792855890229) -- (-2.212179279749624,0.559879936472536) -- (-2.2121788605924286,0.5598812382383875) -- (-2.2121780222706695,0.5598838417653925) -- (-2.212176345597673,0.5598890488006096) -- (-2.2121729921337705,0.5598994627958729) -- (-2.2121662847343306,0.5599202904857162) -- (-2.212152868048953,0.5599619446626724) -- (-2.212126027132531,0.560045248205688) -- (-2.21207231511962,0.5602118360483308) -- (-2.2119647703943155,0.5605449347616497) -- (-2.2117491983120727,0.5612108243131146) -- (-2.211316124951419,0.5625413720169141) -- (-2.210442272088383,0.5651975426404694) -- (-2.208663826938933,0.5704901912512325) -- (-2.2049846600948473,0.5809967699259874) -- (-2.197142669412508,0.6016954690773643) -- (-2.1795676535637605,0.641838362499196) -- (-2.1488852271622356,0.6985667221527735) -- (-2.1131182778754196,0.7516202447141427) -- (-2.072596494315117,0.801024098493835) -- (-2.0276495650931325,0.8468034518023823) -- (-1.9786071788212705,0.8889834729503155) -- (-1.9257990241113352,0.9275893302481666) -- (-1.8695547895751312,0.9626461920064666) -- (-1.810204163824463,0.9941792265357473) -- (-1.748076835471135,1.02221360214654) -- (-1.6835024931269515,1.0467744871493765) -- (-1.6168108254037172,1.0678870498547877) -- (-1.5483315209132362,1.0855764585733054) -- (-1.4783942682673135,1.099867881615461) -- (-1.4073287560777528,1.110786487291786) -- (-1.335464672956359,1.1183574439128123) -- (-1.2631317075149364,1.1226059197890705) -- (-1.1906595483652898,1.1235570832310928) -- (-1.118377884119223,1.1212361025494102) -- (-1.0466164033885412,1.1156681460545546) -- (-0.9757047947850481,1.1068783820570571) -- (-0.9059725646466279,1.0948920169596625) -- (-0.8376779493738902,1.0797491513951294) -- (-0.7709000070282599,1.0615273312851934) -- (-0.7056900828866213,1.0403098940620525) -- (-0.6420995222258608,1.016180177157904) -- (-0.5801796703228628,0.9892215180049442) -- (-0.5199818724545127,0.9595172540353702) -- (-0.4615574738976953,0.9271507226813795) -- (-0.404957819929296,0.892205261375169) -- (-0.3502342558261997,0.8547642075489355) -- (-0.29743812686529125,0.8149108986348768) -- (-0.24662077832345686,0.7727286720651889) -- (-0.19783355547758008,0.7283008652720697) -- (-0.15112780360454758,0.6817108156877165) -- (-0.1065548679812438,0.6330418607443253) -- (-0.06416609388455363,0.5823773378740944) -- (-0.024012826591361947,0.5298005845092202) -- (0.01385358862144459,0.47539493808190003) -- (0.04938180647698287,0.4192437360243304) -- (0.08252048169836623,0.3614303157687091) -- (0.11321826900871024,0.30203801474723213) -- (0.14142382313112956,0.24115017039209818) -- (0.1670857987887402,0.17885012013550283) -- (0.19015285070465504,0.11522120140964431) -- (0.2105736336019901,0.05034675164671909) -- (0.22829680220386006,-0.015689891721076377) -- (0.24327101123338002,-0.08280539126154296) -- (0.25544868049466674,-0.15090533755381808) -- (0.26481277543972537,-0.21980549546892236) -- (0.27136109372230566,-0.2892780127632677) -- (0.2750915312208053,-0.3590947483436424) -- (0.2760019838136203,-0.4290275611168508) -- (0.2740903473791456,-0.498848309989699) -- (0.2693545177957808,-0.5683288538689735) -- (0.26179239094192,-0.6372410516614786) -- (0.2514018626959582,-0.7053567622740204) -- (0.23818082893629677,-0.7724478446133873) -- (0.22212718554132627,-0.8382861575863907) -- (0.20323882838944796,-0.9026435600998206) -- (0.18151365335905778,-0.9652919110604778) -- (0.15694955632854812,-1.0260030693751667) -- (0.12954443317632025,-1.0845488939506844) -- (0.09929617978076832,-1.14070124369383) -- (0.0662026920202865,-1.1942319775114059) -- (0.03026186577327694,-1.2449129543102004) -- (-0.008528403081868863,-1.2925160329970304) -- (-0.05017021866675231,-1.3368130724786753) -- (-0.09466568510297657,-1.3775759316619536) -- (-0.14201642639025103,-1.4145771986379412) -- (-0.1921612955629559,-1.4476847947890281) -- (-0.24491463536845437,-1.4769557412734091) -- (-0.3000761134865275,-1.502469346981977) -- (-0.35744539759689964,-1.5243049208056032) -- (-0.41682215537935186,-1.5425417716351628) -- (-0.47800605451362976,-1.557259208361545) -- (-0.5407967626795145,-1.5685365398756108) -- (-0.6049939475567445,-1.57645307506826) -- (-0.6703972768250726,-1.5810881228303604) -- (-0.7368064181642726,-1.5825209920527872) -- (-0.8040210392540903,-1.5808309916264296) -- (-0.8718408077742996,-1.576097430442159) -- (-0.940065391404632,-1.5683996173908614) -- (-1.0084944578248827,-1.557816861363408) -- (-1.0769276747147618,-1.5444284712506704) -- (-1.145164709754079,-1.528313755943552) -- (-1.2130052306225512,-1.5095520243329208) -- (-1.280248904999974,-1.4882225853096305) -- (-1.3467049008445136,-1.4644048101305351) -- (-1.4122318721796097,-1.4381783949107891) -- (-1.476704534147756,-1.4096231412010773) -- (-1.5399976103111461,-1.378818850607363) -- (-1.6019858242319875,-1.345845324735592) -- (-1.662543899472456,-1.310782365191745) -- (-1.7215465595947412,-1.2737097735817748) -- (-1.7788685281610617,-1.2347073515116271) -- (-1.8343845287336,-1.1938549005872794) -- (-1.887969284874532,-1.1512322224146807) -- (-1.9394975201460651,-1.1069191185997944) -- (-1.9888439581104072,-1.0609953907485625) -- (-2.0358833223297337,-1.0135408404669626) -- (-2.080490336366238,-0.9646352693609437) -- (-2.1225397237821246,-0.9143584790364834) -- (-2.161906208139576,-0.8627902710995095) -- (-2.1984645130007996,-0.8100104471559995) -- (-2.232089361927968,-0.7560988088119096) -- (-2.2626554784833024,-0.7011351576732174) -- (-2.2900375862289755,-0.6451992953458436) -- (-2.3141104087271804,-0.5883710234357729) -- (-2.334748669540132,-0.5307301435489542) -- (-2.3518270922300095,-0.4723564572913581) -- (-2.365220400358985,-0.4133297662689408) -- (-2.374804412420559,-0.35372974379127964) -- (-2.3805786720811284,-0.2936215659212271) -- (-2.3827784417801467,-0.23304278885132046) -- (-2.3816654389257224,-0.1720278689667003) -- (-2.3775013809257928,-0.11061126265250376) -- (-2.3705479851884377,-0.04882742629388925) -- (-2.361066969121623,0.013289183724016773) -- (-2.3493200501334854,0.075704111016063) -- (-2.335568945631934,0.13838289919710878) -- (-2.32007537302502,0.20129109188200278) -- (-2.3031010497207944,0.264394232685615) -- (-2.2849076931272236,0.32765786522279416) -- (-2.265757020652444,0.39104753310839957) -- (-2.2459107497043647,0.45452877995728347) -- (-2.2256305976910085,0.5180671493843159) -- (-2.2154096029749155,0.54984698744427) -- (-2.212851462904098,0.5577923276252195) -- (-2.2122118869465908,0.5597786680603862) -- (-2.212191900149918,0.5598407412053072) -- (-2.212181906751482,0.5598717777777802) -- (-2.2121806575766527,0.5598756573493384) -- (-2.212180032989295,0.5598775971351175) -- (-2.2121797206955307,0.5598785670280089)(0.20710029832073262,0.06209629062214206);
\draw[line width=1.pt,color=grey1,fill=grey1,fill opacity=0.1] (4.2839013914239485,2.5886583736176347) -- (4.294644748608517,2.5930310566578445) -- (4.364942593212851,2.6181840380144448) -- (4.4336322658900835,2.636876348115444) -- (4.566557784254766,2.6561780237328243) -- (4.694162681285618,2.653534221874356) -- (4.817188334565688,2.6315430809044096) -- (4.936376121678026,2.5928027391873565) -- (5.052467420205684,2.539911335087568) -- (5.1662036077317115,2.4754670069694145) -- (5.278326061839159,2.402067893197267) -- (5.389576160111078,2.322312132135497) -- (5.500695280130517,2.238797862148475) -- (5.612424799480529,2.154123221600572) -- (5.725506095744163,2.0708863488561597) -- (5.840680546504469,1.991685382279608) -- (5.9586895293444995,1.919118460235288) -- (6.080274421847303,1.8557837210875712) -- (6.2061766015959305,1.8042793032008282) -- (6.336919862556265,1.7666586363467136) -- (6.471832083494934,1.7419812425124244) -- (6.540513814910263,1.733882680782206) -- (6.60983242861754,1.7282834493093313) -- (6.679652775862731,1.7249379781933465) -- (6.749839707891802,1.7236006975337927) -- (6.8202580759507185,1.7240260374302192) -- (6.890772731285443,1.725968427982166) -- (6.961248525141944,1.7291822992891808) -- (7.031550308766184,1.733422081450807) -- (7.101542933404131,1.7384422045665895) -- (7.171091250301746,1.7439970987360747) -- (7.308314365859853,1.7557289206343212) -- (7.442138465408222,1.7666529879419057) -- (7.571482358914577,1.7748047414551875) -- (7.695264856346639,1.7782196219705204) -- (7.812404767672128,1.77493307028426) -- (7.921820902858768,1.7629805271927745) -- (8.022432071874277,1.7403974334924035) -- (8.1131669347531,1.7052405078154607) -- (8.193452162796849,1.6566422586677305) -- (8.263437151330375,1.5952964025774783) -- (8.323328051352064,1.5220192580759075) -- (8.373331013860295,1.4376271436942254) -- (8.413652189853448,1.342936377963639) -- (8.444497730329914,1.2387632794153554) -- (8.466073786288067,1.1259241665805853) -- (8.478586508726298,1.005235357990525) -- (8.482242048642984,0.8775131721763891) -- (8.477246557036512,0.743573927669388) -- (8.471569096599087,0.6745280079470604) -- (8.463806184905263,0.6042339430007182) -- (8.453983590829843,0.5327937726467624) -- (8.442127083247616,0.46030953670159747) -- (8.428262431033392,0.38688327498161) -- (8.412415403061967,0.31261702730321517) -- (8.39461176820813,0.2376128334828138) -- (8.374877295346682,0.16197273333679263) -- (8.353237753352424,0.08579876668157382) -- (8.329718911100159,0.009192973333536969) -- (8.304346537464667,-0.06774260689089928) -- (8.277146401320767,-0.1449059341753447) -- (8.24814427154324,-0.22219496870339484) -- (8.217365917006894,-0.29950767065865946) -- (8.184837106586528,-0.3767420002247128) -- (8.150583609156932,-0.4537959175851789) -- (8.114631193592906,-0.5305673829236461) -- (8.077005628769253,-0.6069543564237172) -- (8.037732683560769,-0.6828547982689948) -- (7.996838126842243,-0.7581666686430601) -- (7.954347727488483,-0.8327879277295303) -- (7.910287254374286,-0.9066165357119935) -- (7.864682476374445,-0.9795504527740597) -- (7.817559162363769,-1.0514876390993315) -- (7.7689430812170395,-1.1223260548713831) -- (7.718860001809067,-1.1919636602738315) -- (7.614395923708561,-1.3272282807043219) -- (7.504373079060646,-1.4564651818595742) -- (7.388997618863692,-1.5788580452084098) -- (7.268475740362682,-1.6935906402034426) -- (7.143102321562111,-1.8000154512002524) -- (7.013445608776898,-1.8980050450167596) -- (6.880126355700931,-1.9875318836158868) -- (6.743765316028082,-2.0685684289605533) -- (6.674638111747136,-2.1058942360258754) -- (6.604983243452189,-2.141087143013685) -- (6.534878305354958,-2.1741437076693373) -- (6.464400891667147,-2.2050604877382014) -- (6.393628596600543,-2.233834040965643) -- (6.322639014366857,-2.2604609250970285) -- (6.251509739177784,-2.284937697877723) -- (6.180318365245128,-2.307260917053089) -- (6.109142486780563,-2.3274271403684956) -- (6.038059697995877,-2.3454329255693054) -- (5.967147593102773,-2.3612748304008875) -- (5.896483766312947,-2.3749494126086006) -- (5.826145811838238,-2.3864532299378176) -- (5.756211323890263,-2.395782840133897) -- (5.686757896680838,-2.4029348009422087) -- (5.6178631244216675,-2.407905670108118) -- (5.4820599216009995,-2.411290364494185) -- (5.349422469122182,-2.4059093852550255) -- (5.220571520679066,-2.391735194353558) -- (5.0961278299655035,-2.3687402537527014) -- (4.976712150675404,-2.3368970254153894) -- (4.862945236502632,-2.2961779713045374) -- (4.7554478411410415,-2.2465555533830717) -- (4.654840718284525,-2.188002233613915) -- (4.561744621626907,-2.1204904739599826) -- (4.47677982832829,-2.0439929745392362) -- (4.4001152707214715,-1.9587080019685175) -- (4.330629435250273,-1.865478742986454) -- (4.26697166212216,-1.7652629036772964) -- (4.207791291544652,-1.659018190125039) -- (4.151737663725328,-1.5477023084140171) -- (4.097460118871766,-1.4322729646283108) -- (4.043607997191373,-1.3136878648520849) -- (3.988830638891727,-1.1929047151695897) -- (3.931777384180407,-1.0708812216648766) -- (3.8710975732648762,-0.9485750904222243) -- (3.8054405463526564,-0.8269440275257409) -- (3.7334556436513253,-0.7069457390596199) -- (3.653794357893389,-0.5895371128871076) -- (3.5658963526686875,-0.47537543621387357) -- (3.4711602604954805,-0.3643733501150379) -- (3.3712839139101334,-0.2563297633178294) -- (3.2679651454496934,-0.15104358454944844) -- (3.162901787650412,-0.048313722537095316) -- (3.0577916730487686,0.0520609139920154) -- (2.9543326341812417,0.15028141631069758) -- (2.8542225035846513,0.24654887569170825) -- (2.759159113795249,0.34106438340786127) -- (2.6708402973497414,0.43402903073199894) -- (2.590963886784266,0.5256439089368357) -- (2.521227714635984,0.6161101092951995) -- (2.463329613441033,0.7056287230799185) -- (2.418967415736006,0.7944008415637498) -- (2.389838954057268,0.8826275560194645) -- (2.377540942512269,0.9705028113993137) -- (2.382228407342268,1.0581186645647112) -- (2.402981953556832,1.1454912400407657) -- (2.438856517544764,1.232634848279794) -- (2.4889070356958882,1.3195637997341105) -- (2.552188444398894,1.4062924048560266) -- (2.627755680042924,1.4928349740978621) -- (2.714663679017235,1.5792058179119284) -- (2.8119673777108574,1.6654192467505347) -- (2.918721712513161,1.751489571066001) -- (3.0339816198131757,1.8374311013106386) -- (3.156802035999931,1.923258147936762) -- (3.2862378974626836,2.008985021396688) -- (3.421344140590463,2.0946260321427292) -- (3.490728322900509,2.13741881101439) -- (3.5611757017726404,2.1801954906272005) -- (3.632568144255629,2.222957859787696) -- (3.7047875173982447,2.265707707302413) -- (3.7777156882488043,2.3084468219778973) -- (3.851234523856192,2.351176992620685) -- (3.925225891269065,2.393900008037315) -- (3.9995716575359666,2.4366176570343274) -- (4.074153689705895,2.479331728418261) -- (4.1114963889566525,2.500687981507367) -- (4.130174198978125,2.5113660102265616) -- (4.139513911538188,2.51670501235799) -- (4.1441838687619565,2.5193745118951814) -- (4.146518859991829,2.5207092614727156) -- (4.2839013914239485,2.5886583736176347) -- (4.283896524216145,2.5886563611255324);

\draw [line width=1.pt,color=fgreen,fill=fgreen,fill opacity=0.1] (5.8,-0.4) circle (1.6cm);
\draw [rotate around={-52.1250163489031:(-0.923941480585115,-0.5762138995567037)},line width=1.pt,color=fgreen,fill=white,fill opacity=1.0] (-0.9,-0.6) ellipse (0.7cm and 0.7cm);
\draw [line width=1.pt,color=fgreen, fill=white] (6.1,-0.6) circle (0.9cm);

\draw (-3.4,3.0) node[anchor=north west] {$U$};
\draw (5.5,2.8) node[anchor=north west] {$V$};
\draw (6.7,1.3) node[anchor=north west] {$B'$};
\draw (6.2,0.2) node[anchor=north west] {$B$};
\draw (-2.3,1.5) node[anchor=north west] {$\Lambda'$};
\draw (-1.4,0.1) node[anchor=north west] {$\Lambda$};

\draw [->,line width=0.5pt] (-0.6383519974985363,1.5450755326920593) to [bend left = 30] node[above]{at most $D$ singular values} (2.7385498928069865,1.6012261499452407);
\draw [->,line width=0.5pt] (0.15407289126178503,0.21139129699848636) to [bend left = 20] node[above, xshift=-11pt]{$d:1$} (4.4,0.38);
\draw [->,line width=0.5pt] (-0.19259806846984373,-0.6862772943356041) to [bend left = 10] (5.20843113253199,-0.6543005980134807);
\end{tikzpicture}
\end{center}
\caption{The domains in the Covering Lemma.}
\label{Fig:Covering}
\end{figure}
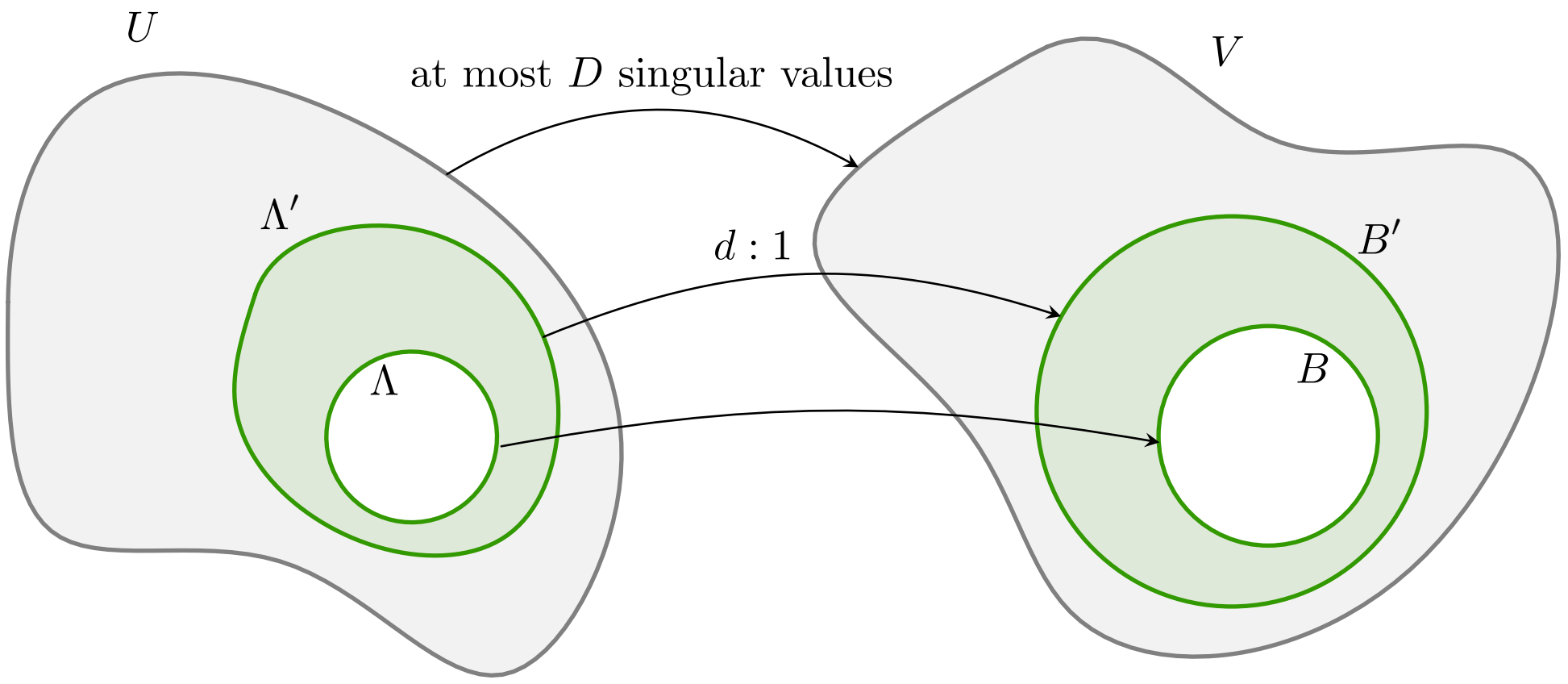

The Covering Lemma shows that if the annuli are sufficiently degenerate (see Figure~\ref{Fig:Covering}), then their pullback is essentially controlled by the local degree $d$, which may be \emph{much smaller} than the number of sheets over $V \sm S_f$; the latter may be even infinite. This theorem recovers the statement of \cite{KL} when $f \colon U \to V$ has finite degree (up to a universal multiplicative constant). Let us give some further remarks.

\begin{remark}[Value of $\eps$]
    It will be clear from the proof that we can choose 
    \[
    \eps = \frac{\eta}{32\pi D^3}.
    \]
\end{remark}

\begin{remark}
As in \cite{KvS, CDKS}, the Covering Lemma can be stated as
\[
\modulus\left(U \sm \ovl \Lambda\right) \ge \min \left\{\eps, \frac{1}{\eta} \modulus\left(B' \sm \ovl B\right), \frac{\eta}{32d^2} \modulus\left(V \sm \ovl B\right)\right\}.
\]
\end{remark}

\begin{remark}[Finite degree Covering Lemma]
\label{Rem:Rational}
If $f \colon U \to V$ is proper of degree $D$, then our approach gives that, provided
\[
\modulus\left(U \sm \ovl \Lambda\right) < \frac{1}{32\pi D^3} =: \eps,
\]
we have the following estimate:
\[
\modulus \left(V \sm \ovl B\right) \modulus (B' \sm \ovl B) < 32d^2 \modulus\left(U \sm \ovl \Lambda\right)^2.
\]
This gives the desired bound under the collar assumption $\cmod{B' \sm \ovl B} \ge \eta\cmod{U \sm \ovl \Lambda}$. In particular, the degeneration parameter $\eps$ does not depend on the collar parameter $\eta$. 
\end{remark}

\begin{remark}
    The proof presented in this paper is motivated by the recent successful applications of the hyperbolic area method in the new proof of the absence of wandering domains in rational dynamics (Sullivan's No Wandering Domains Theorem \cite{Sul85}, see \cite{ye26}) and the absence of bounded orbit wandering domains for entire, meromorphic, and, more crucially, locally defined holomorphic maps \cite{ye26-2, DPU, PRW}.
\end{remark}

\begin{remark}[Further generalizations]
    We believe that the argument in this paper may be further extended to establish the Quasi-Additivity and Quasi-Invariance Laws of \cite{KL} in terms of extremal width.
\end{remark}

\medskip
\noindent
{\bf AI assistance.} ChatGPT was used in exploratory discussions regarding the content of this manuscript based on our previous work \cite{DPU}. The final results and text were written by the author, who takes all responsibility (and blame, if necessary) for the content. 

\section{Some preliminary lemmas}

Let $X$ be a hyperbolic Riemann surface. We will denote by $\rho_X$ the local density of the complete hyperbolic metric on $X$ (normalized to have curvature $-1$). For a curve $\gamma \subset X$, we will denote 
\[
\ell_X(\gamma) = \int_\gamma \rho_X(z) |dz|
\]
the hyperbolic length of this curve; and similarly, for a measurable subset $E \subseteq X$, we will write
\[
\Area_X(E) = \int_E \rho_X^2(z) |dz|^2
\]
its hyperbolic area; here, $z=x+iy$ is the local holomorphic coordinate on $X$, and $|dz|^2 = dxdy$ denotes the Euclidean area element.

We will need the concept of Dirichlet integrals and their relation to the conformal modulus of an annulus; see~\cite[Appendix]{KL}. See also \cite[\S6.4]{LyubichBook} for the introduction to this concept. We specify the definitions suitably for our purposes. Let $Y \Subset X \subset \C$ be a pair of open, compactly contained Jordan disks. Consider a continuous function $h \colon \ovl X \to \R$ that is equal to $1$ on $\ovl Y$, equal to $0$ on $\partial X$, and is harmonic in $X \sm \ovl Y$; we will call $h$ the \emph{equilibrium potential} of the pair $(X, Y)$. For such a potential, we can consider its \emph{Dirichlet energy} integral
\[
\int_X |\nabla h(z)|^2\, |dz|^2 = \int_{X \sm \ovl Y} |\nabla h(z)|^2\, |dz|^2.
\]
It is easy to see that this integral does not depend on the choice of a conformal metric on $X$, i.e., if $\rho(z) |dz|$ is a new smooth metric on $X$, then
\begin{equation}
\label{Eq:Conf}
\int_X |\nabla h(z)|^2\, |dz|^2 = \int_X \|\nabla_\rho h(z)\|^2_\rho \cdot \rho(z)^2\, |dz|^2,
\end{equation}
where $\|\cdot\|_\rho$ and $\nabla_\rho$ denote the norm and gradient in the metric $\rho$. And therefore, it is also invariant under conformal transformations of $X$. For $t \in (0,1)$, we will write $\{h=t\}$ for the level curve in $X \sm \ovl Y$ corresponding to $h=t$. This curve is smooth. We adopt a similar notation for the super- and sub-level sets of $h$. 

Since $\Delta h = 0$ in $X \sm \ovl Y$, the divergence theorem on the set $\{a < h < b\} \subset X \sm \ovl Y$ gives
\[
0 = \int_{a < h < b} \Delta h(z) \, |dz|^2 = \int_{\{h = b\}} |\nabla h(z)| \, |dz| - \int_{\{h = a\}} |\nabla h(z)| \, |dz|. 
\]
Thus, the value of the integral $\int_{\{h = t\}} |\nabla h(z)| |dz|$ does not depend on the choice of $t$. Therefore, by the co-area formula,
\begin{equation}
\label{Eq:Indep0}
    \int_X |\nabla h(z)|^2 \, |dz|^2 = \int_0^1 \left(\int_{\{h = t\}}|\nabla h| \, |dz|\right) dt = \int_{\{h = t\}}|\nabla h| \, |dz| \qquad \forall t \in (0,1).
\end{equation}
Note that the last identity is also independent of the choice of a conformal metric on $X$. In particular,
\begin{equation}
    \label{Eq:Indep}
    \int_X |\nabla h(z)|^2 \, |dz|^2 = \frac{1}{b-a} \int_{\{a < h < b\}} |\nabla h(z)|^2 \, |dz|^2 \qquad \forall 0 \le a < b \le 1 
\end{equation}

Finally, the Dirichlet integral of the equilibrium potential of $(X,Y)$ is related to the conformal modulus $\cmod{X \sm \ovl Y}$ of the annulus $X \sm \ovl Y$ as
\begin{equation}
    \label{Eq:Dir}
    \frac{1}{\cmod{X \sm \ovl Y}} = \int_{X} |\nabla h|^2 \, |dz|^2 = \inf \limits_{\substack{v\colon \ovl X \to \R \\ v|_{\ovl Y} =1}} \int_X |\nabla v(z)|^2 \, |dz|^2,     
\end{equation}
where the last infimum is taken over all locally Lipschitz functions $v$ with compact support in $X$ that are equal to $1$ on $\ovl Y$. The last equality comes from the variational viewpoint on the Dirichlet problem. We refer to \cite[Appendix]{KL} and the classical source \cite{ahlfors} for the background on conformal moduli.

\begin{lemma}[\cite{DPU}]
    \label{Lem:GB}
    Let $X \subset \C$ be a topological disk, and $Q \subset X$ be a finite collection of points. Then for every measurable set $E \subseteq X \sm Q$
    \begin{equation}
        \label{Eq:GB}
        0 \le (\Area_{X \sm Q} - \Area_X)(E) \le 2\pi \cdot \#Q.
    \end{equation}
    Here, the area gain is understood in the integral form 
    \[
    (\Area_{X \sm Q} - \Area_X)(E) := \int_E \left(\rho_{X \sm Q}^2(z) - \rho_X^2(z)\right) \, |dz|^2. 
    \]
\qed
\end{lemma}

\begin{lemma}
    \label{Lem:Lenth}
    Let $Y \Subset X$ be a pair of open Jordan disks, and let $\gamma \subset X \sm \ovl Y$ be a smooth simple closed curve enclosing $Y$. Then
    \begin{equation}
        \label{Eq:Length}
        \frac{1}{\cmod{X \sm \ovl Y}} \le \ell_X(\gamma) + \frac{2\pi}{3}.
    \end{equation}
\end{lemma}

\begin{proof}
Since the modulus and the hyperbolic length are conformal invariants, we can assume that $X$ is the hyperbolic plane (modeled in the upper half plane $\H$ with metric $\rho_\H(z) = 1/y$, $z = x + iy$). 

Let $D \Subset \H$ be a domain with rectifiable boundary. Then  
\begin{equation}
\label{Eq:Ineq1}
\Area_{\H}(D) \le \ell_{\H}(\partial D).    
\end{equation}
Indeed, since $D$ can be approximated by domains with smooth boundaries so that the length and area converge, it is enough to show \eqref{Eq:Ineq1} for domains with smooth boundaries $\partial D$. Let $(\partial_x, \partial_y)$ denote the standard basis in each tangent plane of $\H$, then a direct computation shows that the vector field $-y \cdot \partial_y$ has hyperbolic length $1$ and divergence $1$. Let $\nu$ be the outward-facing normal vector field along $\partial D$. Set $dA_\H := \rho_\H^2 \, |dz|^2$ and $ds_\H := \rho_\H \, |dz|$ to be the hyperbolic area and length element in $\H$, and $\left<\cdot, \cdot\right>_\H$, $\|\cdot\|_\H$ be the hyperbolic dot product and length. Using the divergence theorem,
\[
\Area_\H(D) = \int_D \diver(-y\,\partial_y) \,dA_\H = \int_{\partial D} \left<\nu, -y \partial_y\right>_\H \, ds_\H \le \int_{\partial D} \|-y \partial_y\|_\H ds_\H = \ell_\H(\partial D). 
\]
\begin{remark}
    Inequality \eqref{Eq:Ineq1} also follows from the isoperimetric inequality in the hyperbolic plane that has the form $\ell_\H^2(\partial D) \ge 4\pi \Area_\H(D) + \Area_\H^2(D)$.
\end{remark}

Let $\Gamma$ be the convex hull of $\gamma$ (i.e., the smallest hyperbolically convex set in $\H$ that contains $\gamma$). By definition, $\Gamma$ is compact, $\ovl Y \subset \Gamma$. Since $\gamma$ is smooth, $\partial \Gamma$ is rectifiable and $\ell_\H(\partial \Gamma) \le \ell_\H(\gamma)$. In particular, using \eqref{Eq:Ineq1},
\begin{equation}
    \label{Eq:Ineq2}
    \Area_\H(\Gamma) \le \ell_\H(\partial \Gamma) \le \ell_\H(\gamma).
\end{equation}

Let $\Gamma_s$ be the outer parallel set for $\Gamma$ at distance $s \ge 0$, i.e., $\Gamma_s := \{z \in \H : \dist_\H(z, \Gamma) \le s\}$. By construction, $\Gamma = \Gamma_0$. By the Steiner formula in the hyperbolic plane (see e.g.\cite[Chapter 18]{santalo2004integral}),
\begin{equation}
    \label{Eq:Stainer}
    \ell_\H(\partial \Gamma_s) = \ell_\H(\partial \Gamma) \cosh s + \left(2\pi + \Area_\H(\Gamma)\right) \sinh s.
\end{equation}

Consider the function $v(z) := e^{-{\dist_\H(z,\Gamma)}}$. It satisfies $v|_{\Gamma} = 1$ and $v(z)\to 0$ when $z \to \partial\H$. Since $v$ is not compactly supported in $\H$, for $\eps \in (0,1)$, we define its cutoff $v_\eps(z) := \max\{v(z) - \eps, 0\} / (1-\eps)$. This function is locally Lipschitz, $v_\eps|_\Gamma = 1$, and the support of $v_\eps$ is precisely equal to $\Gamma_{\log(1/\eps)}$.

By conformal invariance \eqref{Eq:Conf} and the variational approach \eqref{Eq:Dir},
\[
\frac{1}{\cmod{\H \sm \ovl Y}} \le \int_{\H} \|\nabla_{\rho_\H} v_\eps(z)\|_{\rho_\H}^2 dA_\H = \frac{1}{(1-\eps)^2} \int_{\{v > \eps\}} \|\nabla_{\rho_\H} v(z)\|_{\rho_\H}^2 dA_\H.
\]
The right-hand side converges, as $\eps \to 0$, to $\int_{\H} \|\nabla_{\rho_\H} v(z)\|_{\rho_\H}^2 dA_\H$, thus giving the estimate
\begin{equation}
    \label{Eq:Dir2}
    \frac{1}{\cmod{\H \sm \ovl Y}} \le \int_{\H} \|\nabla_{\rho_\H} v(z)\|_{\rho_\H}^2 dA_\H    
\end{equation}
Note that $\nabla_{\rho_\H} v(z)=-v(z) \cdot \nabla_{\rho_\H} \dist_\H(z,\Gamma)$. Since the gradient of a hyperbolic distance function to a convex set always has hyperbolic norm one outside of that set, we obtain that $\|\nabla_{\rho_\H} v(z)\|_{\rho_\H} = v(z) = e^{-\dist_\H(z, \Gamma)}$ everywhere outside of $\Gamma$. By the co-area formula,
\begin{equation}
    \label{Eq:Est-1}
    \int_{\H} \|\nabla_{\rho_\H} v(z)\|_{\rho_\H}^2 dA_\H = \int_0^1  \int_{\{v(z) = t\}} \|\nabla_{\rho_\H} v(z)\|_{\rho_\H} \, ds_\H dt   
\end{equation}
We can now compute the right-hand side using the Steiner formula \eqref{Eq:Stainer}:
\begin{equation*}
    \begin{aligned}
        \int_0^1  \int_{\{v(z) = t\}} \|\nabla_{\rho_\H} v(z)\|_{\rho_\H} \, ds_\H dt &= \int_0^1 t \cdot \ell_\H(\{v(z)=t\}) dt = \int_0^1 t \cdot \ell_\H(\partial\Gamma_{-\log t}) dt\\
        & =\frac{1}{2}\int_0^1 \left((1+t^2)\ell_\H(\partial \Gamma) + (1-t^2)(2\pi + \Area_\H(\Gamma))\right)dt \\
        &= \frac{2}{3}\ell_\H(\partial \Gamma) + \frac{1}{3}(2\pi + \Area_\H(\Gamma)) \le \ell_\H(\gamma) + \frac{2\pi}{3},
    \end{aligned}    
\end{equation*}
where the last inequality follows by \eqref{Eq:Ineq2}. Combining the last estimate with \eqref{Eq:Dir2}, we obtain $\left(\cmod{\H \sm \ovl Y}\right)^{-1} \le \ell_\H(\gamma) + 2\pi/3$, as required.
\end{proof}

\begin{lemma}
    \label{Lem:Area}
    Let $Y \Subset X$ be a pair of open Jordan disks, and let $h$ be the equilibrium potential of $(X,Y)$. Then, for $t \in (0,1)$,
    \begin{equation}
        \label{Eq:Area}
        \Area_X(\{h > t\}) \le \frac{1}{t \cdot \cmod{X \sm \ovl Y}}.       
    \end{equation}
    
\end{lemma}

\begin{proof}
Recall that by \eqref{Eq:Indep0} and \eqref{Eq:Dir}, for every $t \in (0,1)$ we have
\[
\frac{1}{\cmod{X \sm \ovl Y}} = \int_{h=t} |\nabla h(z)| \, |dz|.
\]
Set $a(t) := \Area_X(\{h > t\})$. By the co-area formula,
\begin{equation*}
\begin{aligned}
    a(t) &=  \Area_X(Y) + \int_{\{h > t\}} \rho_X(z)^2 \, |dz|^2 = \int_{t}^1 \int_{h=s} \frac{\rho_X^2(z)}{|\nabla h(z)|} \, |dz| ds \quad \Longrightarrow \\ 
    &\Longrightarrow\quad a'(t) = -\int_{\{h=t\}} \frac{\rho_X^2(z)}{|\nabla h(z)|} \, |dz|.   
\end{aligned}
\end{equation*}
By the Cauchy--Schwarz inequality, 
\begin{equation*}
    \begin{aligned}
        \big(\ell_X(\{h = t\})\big)^2 &= \left(\int_{\{h=t\}}\frac{\rho_X(z)}{\sqrt{|\nabla h(z)|}}\sqrt{|\nabla h(z)|} \, |dz|\right)^2 \\
        &\le \int_{\{h=t\}} \frac{\rho_X^2(z)}{|\nabla h(z)|} \, |dz| \,\cdot  \int_{\{h=t\}} |\nabla h(z)| \, |dz| = - \frac{a'(t)}{\cmod{X \sm \ovl{Y}}} 
    \end{aligned}
\end{equation*}
On the other hand, since $\partial\{h > t\} = \{h = t\}$, inequality \eqref{Eq:Ineq1} yields:
\[
a^2(t) = \left(\Area_{X}(\{h > t\})\right)^2 \le \left(\ell_X({\{h = t\}})\right)^2 \le -\frac{a'(t)}{\cmod{X \sm \ovl{Y}}}.
\]
This differential inequality is equivalent to $(1/a(t))' \ge \cmod{X \sm \ovl Y}$. Integrating this inequality from $s$ to $t$, where $0 < s < t$ gives $1/a(t) \ge  1/a(s) + (t-s) \cmod{X \sm \ovl{Y}} \ge (t-s) \cmod{X \sm \ovl{Y}}$. Sending $s \to 0$ gives the required inequality $1/a(t) \ge t \cmod{X \sm \ovl Y}$.
\end{proof}

\section{Proof of Theorem~\ref{Thm:Main}}
We start by extending the geometric Lemma~\ref{Lem:GB} to the dynamical context. This step is similar to \cite[Section 2]{DPU} (and the analogous steps in \cite{PRW, ye26}).

\begin{lemma}
    \label{Lem:Dynmaical}
    In the notation of Theorem~\ref{Thm:Main}, if $E \subseteq \Lambda'$ is a measurable set, then
    \[
    \Area_U(E) \le d \cdot \Area_{V}(f(E)) + 2\pi d D.
    \]
\end{lemma}
\begin{proof}
    Let $Q \subset V$ be the set of singular values of $f \colon U \to V$. Removing the finite set $E\cap f^{-1}(Q)$ from $E$ does not change the area of $E$ and $f(E)$, so we may assume that $f(E)\cap Q=\emptyset$. 
    
    We have that $f \colon U \sm f^{-1}(Q) \to V \sm Q$ is a covering, and hence a local isometry in the hyperbolic metric. Hence, combined with the domain monotonicity $U \sm f^{-1}(Q) \subset U$, in local densities we have the following estimate: 
    \[
    \rho_{U} \le \rho_{U \sm f^{-1}(Q)} = |f'| \cdot \rho_{V \sm Q} \circ f.
    \]
    This way, we obtain the right-hand side of the estimate
    \begin{equation}
        \label{Eq:Est1}
        0 \le \rho_U^2 - |f'|^2 \cdot \rho_V^2 \circ f \le |f'|^2 \cdot (\rho_{V \sm Q}^2 - \rho_V^2) \circ f;    
    \end{equation}
    the left-hand side follows by the Schwarz--Pick lemma applied to the holomorphic map $f \colon U \to V$. By Lemma~\ref{Lem:GB} applied to the set $B'\sm Q$ gives
    \begin{equation}
        \label{Eq:Est1.5}
        \begin{aligned}
        &0 \le \int_{f(E)} \left(\rho^2_{V \sm Q}(z) - \rho^2_V(z)\right) \, |dz|^2 \le 2\pi \cdot \#Q  \le  2\pi  D \qquad \Longrightarrow \\
        &\qquad \Longrightarrow \qquad\Area_{V\sm Q}(f(E)) \le \Area_V(f(E)) + 2\pi D.    
        \end{aligned}
    \end{equation}
    The multiplicity of $f$ over the set $E$, since $E$ lies in $\Lambda'$, is at most $d$. Therefore, by integrating \eqref{Eq:Est1} over $E$ we obtain:
    \[
    \Area_U(E) \le d\Area_{V} (f(E)) + d \left(\Area_{V\sm Q}(f(E)) - \Area_V(f(E))\right) \le d\Area_{V} (f(E)) + 2\pi d D,
    \]
    where the last inequality follows from \eqref{Eq:Est1.5}. This finishes the proof of the lemma.
\end{proof}

\begin{proof}[Proof of Theorem~\ref{Thm:Main}]
Let $v$ and $b$ be the equilibrium potentials for pairs $(V, B)$ and $(B', B)$ respectively. By the maximum principle, we have that 
\begin{equation}
    \label{Eq:MP}
    v(z) \ge b(z) \quad \forall z \in B'.
\end{equation}
Furthermore, define 
\begin{equation}
    u(z) = \left\{
        \begin{aligned}
        &b \circ f(z), &\quad& z \in \Lambda'\\
        &0, &\quad&z \in U \sm {\Lambda'}
        \end{aligned}
    \right.
\end{equation}
This function has finite Dirichlet energy. 

Define a set $\Omega := \{1/2 < u < 1\} \subset \Lambda' \sm \ovl \Lambda$. Let us apply Lemma~\ref{Lem:Dynmaical} to $\Omega$. Note that $f(\Omega) \subseteq \{b > 1/2\} \subset B'$. By Lemma~\ref{Lem:Dynmaical} followed by inequality \eqref{Eq:MP} and an application of Lemma~\ref{Lem:Area}, we obtain:
\begin{equation}
    \label{Eq:Est2}
    \begin{aligned}
    \Area_U(\Omega) &\le d \,\Area_V(\{b > 1/2\}) + 2\pi d D \le d\,\Area_V(\{v > 1/2\}) + 2\pi d D \\
    &\le \frac{2d}{\cmod{V \sm \ovl B}} + 2\pi d D.
    \end{aligned}
\end{equation}

Observe that 
$
\int_\Omega |\nabla u(z)|^2 \, |dz|^2 = d\int_{\{1/2 < b < 1\}} |\nabla b(z)|^2 \, |dz|^2. 
$
On the other hand, by the constant flux property \eqref{Eq:Indep} and relation \eqref{Eq:Dir}, we have
\[
\int_{\{1/2 < b < 1\}} |\nabla b(z)|^2 \, |dz|^2 = \frac{1}{2}\int_{B'} |\nabla b(z)|^2 \, |dz|^2 = \frac{1}{2\cmod{B' \sm \ovl B}}.
\]
Hence,
\begin{equation}
    \label{Eq:Est3}
    \int_{\Omega} |\nabla u(z)|^2 \, |dz|^2 = \frac{d}{2 \cmod{B' \sm \ovl B}}.
\end{equation}

Finally, for every regular value $t \in (1/2, 1)$ of $u$, the level set $\{u = t\}$ is a finite union of smooth Jordan curves, at least one of which, say $\gamma$, bounds a disk in $U$ that contains $\ovl \Lambda$. By Lemma~\ref{Lem:Lenth} applied to $\Lambda \Subset U$,
\begin{equation}
    \ell_U\left(\{u = t\}\right) \ge \ell_U(\gamma) \ge \frac{1}{\cmod{U \sm \ovl \Lambda}} - \frac{2\pi}{3}.
\end{equation}
Note that there are only finitely many exceptional levels of $u$ as $\{z \in U : 1/2 < u(z) < 1\} \subseteq \Lambda'$ and $f|_{\Lambda'}$ has finite degree. Therefore, by the co-area formula and the Cauchy--Schwarz inequality,
\begin{equation}
    \label{Eq:Est4}
    \begin{aligned}
    \frac{1}{2}\left(\frac{1}{\cmod{U \sm \ovl \Lambda}} - \frac{2\pi}{3}\right) &\le \int_{1/2}^1 \ell_U(\{u = t\}) dt = \int_{1/2}^1 \int_{\{u = t\}} \rho_U(z) \, |dz| dt \\
    &= \int_{\Omega} \rho_U(z) |\nabla u(z)| \, |dz|^2 \le \left(\Area_U(\Omega) \cdot \int_{\Omega} |\nabla u(z)|^2 \, |dz|^2\right)^{\frac{1}{2}}      
    \end{aligned}
\end{equation}
Suppose that the modulus $U \sm \ovl\Lambda$ is degenerate so that
\[
\frac{1}{\cmod{U \sm \ovl \Lambda}} - \frac{2\pi}{3} > \frac{1}{2} \cdot \frac{1}{\cmod{U \sm \ovl \Lambda}} \qquad \Longleftrightarrow \qquad \cmod{U \sm \ovl \Lambda} < \frac{3}{4\pi}.
\]
Then the left-hand side of \eqref{Eq:Est4} is positive, and thus by squaring \eqref{Eq:Est4} and combining with estimates \eqref{Eq:Est2} and \eqref{Eq:Est3}, we obtain:
\begin{equation}
    \label{Eq:Est5}
    \frac{1}{\cmod{U \sm \ovl \Lambda}^2} \le \frac{16d^2}{\cmod{B' \sm \ovl B}} \cdot \left(\frac{1}{\cmod{V \sm \ovl B}} + \pi D\right)
\end{equation}

In Theorem~\ref{Thm:Main}, it is enough to establish the second case. Under the collar assumption $\cmod{B' \sm \ovl B} \ge \eta \cmod{U \sm \ovl \Lambda}$, the estimate \eqref{Eq:Est5} implies
\begin{equation}
\label{Eq:Est6}
\frac{1}{\cmod{U \sm \ovl \Lambda}} \le \frac{16d^2}{\eta \cmod{V \sm \ovl B}} + \frac{16\pi d^2 D}{\eta}.
\end{equation}
Now if $U \sm \ovl \Lambda$ is further degenerate so that 
\[
\frac{1}{\cmod{U \sm \ovl \Lambda}} > \frac{1}{2\cmod{U \sm \ovl \Lambda}} + \frac{16\pi d^2D}{\eta} \qquad \Longleftrightarrow \qquad \cmod{U \sm \ovl \Lambda} < \frac{\eta}{32\pi d^2 D},
\]
then from \eqref{Eq:Est6} we obtain the required bound
\[
\cmod{U \sm \ovl \Lambda} > \frac{\eta}{32d^2} \cmod{V \sm \ovl B}.
\]
This finishes the proof of the covering lemma. Our estimates start to work once the annulus $U \sm \ovl \Lambda$ is $\eps$-degenerate, where we can pick
\[
\eps := \frac{\eta}{32\pi D^3} \le \frac{\eta}{32\pi d^2 D} < \frac{3}{4\pi}.
\]
\end{proof}

\begin{remark}
    If in the statement of Theorem~\ref{Thm:Main} the map $f \colon U \to V$ has a finite degree bounded by $D$, then in the proof above we use the crude bound $\cmod{U \sm \ovl \Lambda} \ge \frac{1}{D} \cmod{V \sm \ovl B}$ (see, e.g., \cite[Lemma 4.5]{KL}) to pass the degeneration from $U \sm \ovl \Lambda$ to $V \sm \ovl B$ and this way bound the second bracket in \eqref{Eq:Est5} without the collar assumption. This gives the statement outlined in Remark~\ref{Rem:Rational}. 
\end{remark}

\bibliographystyle{amsalpha}
\bibliography{references}

\end{document}